\documentclass[12pt,reqno]{amsart}
\usepackage{amscd,amssymb,euscript,enumerate}

\parskip=\baselineskip
\newcommand{\rar}[2]{\xrightarrow[#2]{#1}}
\newcommand{\x}{\times}

\newcommand{\cp}{\rtimes}

\newcommand{\iso}{\cong}

\newcommand{\cc}{\colon}
\newcommand{\co}{\circ}
\newcommand{\oab}{\bar O_A}
\newcommand{\oa}{O_A}
\newcommand{\oat}{O_{A^T}}
\newcommand{\D}{{\raise 0.09em \hbox{/}} \kern -.58em {\partial}}
\newcommand{\dd}{\EuScript{D}}

\newcommand{\te}{\otimes}
\newcommand{\al}{\alpha}
\newcommand{\be}{\beta}
\newcommand{\cs}{C^*}
\newcommand{\ru}{\mathcal{R}^u}
\newcommand{\rs}{\mathcal{R}^s}
\newcommand{\h}{\mathcal{H}}
\newcommand{\B}{\mathcal{B}}
\newcommand{\E}{\mathcal{E}}
\newcommand{\C}{\mathbb{C}}
\newcommand{\R}{\mathbb{R}}
\newcommand{\sk}{\mathcal{S}}
\newcommand{\s}{\sigma}
\newcommand{\tk}{\tau}

\newcommand{\K}{\mathcal{K}}
\newcommand{\N}{\mathcal{N}}
\newcommand{\T}{\mathcal{T}}
\newcommand{\NN}{\mathbb{N}}
\newcommand{\F}{\mathcal{F}}

\newcommand{\de}{\delta_*}
\newcommand{\ds}{\Delta_*}
\newcommand{\Z}{\mathbb{Z}}

\newtheorem{theorem}{Theorem}[section]
\newtheorem{lem}[theorem]{Lemma}
\newtheorem{definition}[theorem]{Definition}
\newtheorem{prop}[theorem]{Proposition}

\theoremstyle{remark}

\begin{document}
\baselineskip=24pt
\title{K-theoretic duality for shifts of finite type}
\author{Jerome Kaminker}
\address{Department of Mathematical Sciences, IUPUI, Indianapolis,
IN 46202-3216} 
\email{jkaminker@math.iupui.edu}
\thanks{The first author was supported in part by NSF Grant DMS-9401457}
\thanks{The second author was supported in part by an NSERC Grant}

\author{Ian Putnam}
\address{Department of Mathematics, University of Victoria, Victoria,
BC V8W-3P4}
\email{putnam@smart.math.uvic.ca}          

 
\begin{abstract}
We will study the stable and unstable Ruelle algebras associated to a
hyperbolic homeomorphism of a compact space.  To do this, we will
describe a notion of K-theoretic duality for $\cs$-algebras which
generalizes Spanier-Whitehead duality in topology.  A criterion for
checking that it holds is presented.  As an application it is shown
that the Ruelle algebras which are associated to subshifts of finite
type, (and agree with Cuntz-Krieger algebras in this case) satisfy
this criterion and hence are dual.
\end{abstract}
\maketitle

\section{Introduction}

This paper is part of a study of certain  $\cs$-algebras which can be associated to a
hyperbolic homeomorphism of a compact space, $(X,f)$.  They are  called  the the stable and unstable
Ruelle algebras, $\rs$ and $\ru$, and are higher dimensional generalizations of
Cuntz-Krieger algebras.  This means that if the dimension of $X$ is
zero, then $\ru \iso \oat \te \K$ and $\rs \iso \oa \te \K$.  One of the basic
results of the theory is a duality relation between $\ru$ and $\rs$.
In the present paper we prove this explicitely in the zero dimensional
case.  Our reason for doing this is to bring out the use of Fock space to construct
the K-theory class implementing the duality in the zero dimensional case.

The notion of Spanier-Whitehead duality in topology has a very natural
generalization to K-theory of $\cs$-algebras.  Briefly, it says that
two algebras, $A$ and $B$, are dual if there are duality classes
$\Delta \in KK^i(A \te B, \C)$ and $\delta \in KK^i(\C, A \te \B)$ which
induce an isomorphism between the K-theory of $A$ and the K-homology of
$B$ via Kasparov product.  It is closely related to the notion of
Poincar\'e duality used by Connes in his study of the standard model
of particle physics, ~\cite{connes:book} .  We will describe it
in more detail in Section 2.  A useful result proved in the paper
is a criterion, presented in Section 3, for deciding when one has a duality
between two algebras.  It is applicable when two duality classes such
as $\Delta$ and $\delta$ are given and one wants to show that they induce 
duality isomorphisms.  Section 4 and Section 5 apply this criterion to
the case of two algebras associated to a hyperbolic dynamical system.
If $A$ is an $n \x n$ aperiodic matrix then one can associate to it
the subshift of finite type, $(\Sigma_A, \s_A)$.  There are two
$\cs$-algebras that can be constructed from this data--the
Cuntz-Krieger algebras $\oa$ and $\oat$.  We show that these algebras
are dual.  In Section 6 we will discuss some further applications and
make some concluding remarks. 

In a later paper we will establish duality for the stable and
unstable Ruelle algebras associated to any hyperbolic homeomorphism of a
compact space, (a Smale space).  Ruelle algebras were introduced by
the second author in ~\cite{ putnam}.  They can be thought of as
higher dimensional generalizations of Cuntz-Krieger algebras.  They
are constructed by defining two equivalence relations on the Smale
space, stable and unstable equivalence.  One takes the $\cs$-algebras
associated to them and then takes the crossed products by the
automorphism induced by the homeomorphism.

The stable and unstable equivalence classes in a Smale space behave
very much like transverse foliations.  Because of that, and the fact
that the homeomorphism is contracting along the stable leaves, one
obtains a duality in K-theory for the algebras.

Cuntz-Krieger algebras are special cases of Ruelle algebras, so the
duality established here would follow from the more general theory.
However, there is an intriguing aspect to this which as yet has no
analogue in the general case. Namely, the duality classes have
representatives constructed using  Fock space.  These classes are obtained
in a natural manner following work of D. Evans, ~\cite{ evans} and
D. Voiculescu, ~\cite{ voiculescu}.  This provides potential
connections with physics (c.f. ~\cite{ jorgensen, dykema-n, faddeev}) and Voiculescu's work on free products which we
hope to pursue in the future.  We would like to thank Dan Voiculescu
and Marius Dadarlat for very helpful conversations.

It should be noted that the general duality theory for Smale spaces
requires a different approach which is based on the notion of
asymptotic morphism.  The final version of the duality theorem uses
these methods, ~\cite{ kaminker-p2}.

\section{K-theoretic duality for $\cs$-algebras}

In this paper we will be describing an example of some $\cs$-algebras
that are dual with respect to K-theory.  This notion of duality has
appeared several times in the past, ~\cite{kasparov:invent,kahn-k-s,parker} and
recently was used by Connes ~\cite{connes:book}.  We present the
definitions here and list some basic facts.  More details can be found
in ~\cite{ kaminker-p1}.  

We will use the following conventions.  Let $\sk$ denote $C_0(\R)$.
Then $KK^1(A,B)$ will be, by definition, equal to $KK(\sk \te A, B)$.
For $A$ and $B$ separable, and $A$ nuclear one has that $KK^1(A,B)
\iso Ext(A,B)$.
We establish some additional notation.  If $\s$ is a permutation, and $A_1, \ldots
A_n$ are algebras, then we will also use $\s$
 to denote the isomorphism
\begin{equation*} 
A_1 \te \cdots \te A_n \to A_{\s (1)} \te \cdots \te A_{\s(n)}.
\end{equation*}  
If $\s$ is a transposition interchanging $i$ and $j$, we will write
$\s_{ij}\cc KK^*( \cdots \te A_i \te \cdots \te A_j \te\cdots , B) \to
KK^*(\cdots \te A_j \te \cdots \te A_i \te \cdots, B) $ for the
homomorphism  induced by $\s$ on the first variable of the Kasparov groups, 
and $\s^{ij}$ for the corresponding map induced on the second variable.  Let $\tk_D \cc KK^i(A,B) \to KK^i(A \te D, B \te D)$ and
$\tk^D \cc KK^i(A,B) \to KK^i(D \te A, D \te B)$ denote the standard
maps, ~\cite{ kasp}.  

Also, we will have need of the following version of Bott
periodicity.  Let 
\begin{equation}
\label{toep} 
\begin{CD}
0 @>>> \K(\ell^2(\NN)) @>>> \T @>{\s_{\T}}>> C(S^1) @>>> 0
\end{CD}
\end{equation}
be the Toeplitz extension.  We will denote it by $\T \in KK^1(C(S^1),
\C) $ and its restriction to $\sk$ by $\T_0 \in KK(\sk \te \sk, \C)$.
Let $\be \in KK(\C, \sk \te \sk)$ be the Bott element.
Then the following holds. (c.f. ~\cite{ blackadar}) 
\begin{theorem} One has
$\be \te_{\sk \te \sk} \T_0 = 1_{\C}$ and $\T_0
\te \be = 1_{\sk \te \sk}$, .
\end{theorem}

We describe next the notion of duality we will be using.
\begin{definition}
Let $A$ and $B$ be  $C^{*}$-algebras.  Suppose that, for $n
= 0$ or $1$, two classes,
$\Delta \in KK^n(A \te B, \C)$ and $\delta \in KK^n(\C, A \te B)$ ,  are given.
Define homomorphisms $\Delta_i \cc K_i(A) \to
K^{i+n}(B)$ and $\delta_i \cc K^{i+n}(B) \to K_i(A)$ in the  following way.
In $n = 1$ set
\begin{equation} 
   \Delta_i(x) = 
\begin{cases}
x \te_{A}\Delta&  \text{if $ i = 0$},\\
\be \te_{\sk \te \sk} ( \s_{12}(x \te_{A} \Delta))&  \text{if $i = 1$} 
\end{cases}
\end{equation}
and let
\begin{equation} 
\delta_i(y) =
\begin{cases}
\be \te_{\sk \te \sk} (\delta \te_{B} y)&  \text{if $i = 0$}, \\
\delta \te_{B} y&  \text{if $ i = 1$},
\end{cases}
\end{equation}

If $n = 0$ set
\begin{equation} 
   \Delta_i(x) = 
\begin{cases}
x \te_{A}\Delta&  \text{if $ i = 0$},\\
 \s_{12}(x \te_{A} \Delta)&  \text{if $i = 1$} 
\end{cases}
\end{equation}
and let
\begin{equation} 
\delta_i(y) =
\begin{cases}
\delta \te_{B} y&  \text{if $i = 0$}, \\
\be \te_{\sk \te \sk} (\delta \te_{B} y)&  \text{if $ i = 1$},
\end{cases}
\end{equation}

We say that $A$ and $B$ are
dual if  
\begin{equation*} 
\Delta_i :K_i(A) \to K^{i+n}(B)
\end{equation*}
and
\begin{equation*} 
\delta_i :K^{i+n}(B) \to K_i(A)
\end{equation*}
are inverse isomorphisms.  Given $A$, if such an algebra $B$ exists it is called a dual of $A$
and it is denoted $\dd A$.    
\end{definition}
In this generality a dual is not unique, so care must taken with the
notation $\dd A$.  We will only use it if a specific dual is in hand.
However, it is easy to see that a dual is unique up to
KK-equivalence.  Indeed, $\s^{12}(\delta ') \te_{A} \Delta \in
KK(B',B)$ and $\s^{12}(\delta) \te_{A} \Delta ' \in KK(B,B')$ yield
the required KK-equivalence.

The form of the definition of the homomorphisms $\Delta_{*}$ and
$\delta_{*}$ is forced by our convention that $KK^1(A,B) = KK(A \te
\sk, B)$.  It is an interesting point that when dealing with an odd
type duality one must bring in some form of Bott periodicity
explicitely.  Either one can incorporate it into the definition of the
homomorphisms as we have done, or one can modify the definitions of
the K-theory groups.  As the reader will see, our choice is the most convenient one for
the proofs we are giving.  Note also that we are working only in the
odd case, (i.e. $n=1$), in this paper.

For a specific algebra $A$ it is not clear
whether a dual, $\dd A$, exists.  In general, the existence of $\dd A$
with prescribed properties, such as separability, is a strong
condition.  If one can take $\dd A$ equal to $A$ then this
agrees with what Connes has developed as Poincar\'e duality in
~\cite{connes:book}.

If one requires only the existence of $\Delta$ and the fact that it
yields an isomorphism in the definition above, then there is no
guarantee that a class $\delta$ exists to give the inverse
isomorphism.  If $A$ was $C(X)$, with $X$ a finite complex, then the
existence of $\delta$ would follow from that of $\Delta$.  However, in
general this need not hold.

The origin of this notion is in Spanier-Whitehead duality in topology,
~\cite{ spanier}.  Recall that if $X$ is a finite complex then there is a dual
complex, $DX$, along with class $\Delta\in H_m(X\wedge
DX)$ 
satisfying that
$\backslash \Delta :H^i(X) \to H_{m-i}(DX)$
is an isomorphism.  The space $DX$ is called the
Spanier-Whitehead dual of $X$.  It is unique up to stable homotopy. If
$M$ is a closed manifold of dimension $n$ embedded in $\R^m$, 
then $DM$ can be taken to be $(\nu M)^{+}$, the Thom space of 
the normal bundle of M. It is interesting to note that there is a relation
between Spanier-Whitehead duality, the Thom isomorphism, $\phi$,  
and Poincar\'e duality 
\begin{equation*} 
\begin{CD}
H_{n-i}(M) @<<{\backslash\Delta}< H^{i+m-n}((\nu M)^+)\\
@A{\cap [M]}AA@AA{\phi}A\\
H^{i}(M)@>>>H^{i}(M),
\end{CD}
\end{equation*}
where $[M] = U \backslash \Delta$, $U$ the Thom class.  Of course,
$\dd (C(X)) = C(D(X))$ for $X$ a finite complex.

If one works in the class $\N$
introduced by Rosenberg and Schochet in their study of the Universal
Coefficient Theorem, ~\cite{ rosenberg-s}, the theory simplifies and
there is a strong analogy with the commutative case.  (However, in general, the
restriction that the algebras lie in $\N$ is too strong.  In several
important examples this does not hold.)  Recall that $ \N$ is defined to be
the smallest class of separable, nuclear $\cs$-algebras  containing $\C$ and
closed under forming extensions, direct limits, and KK-equivalence.
The Universal Coefficient Theorem for KK-theory holds for $KK(A,B)$ if
$B$ is separable and $A \in \N$. Let $\dd\N$ be the subclass of $\N$
consisting of algebras $A$ in $\N$ for which a dual $\dd A$ exists and
is also in $\N$.  For algebras in $\dd\N$ the following facts are easy
consequences of the properties of the Kasparov product and the
Universal Coefficient Theorem.

\begin{enumerate}[i)]
\item If $A$ is dual to $B$, then $B$ is dual to $A$.
\item If $A \in \dd\N$, then  $\dd(\dd A)$ is KK-equivalent to $A$.
\item If $A \in \N$, then $A \in \dd\N$ if and only if $K_*(A)$ is
  finitely generated.
\item Let $E, D \in \N$ and $A \in \dd\N$.  Then
\begin{equation*}
\Delta_* \cc KK^*(E, D \te A) \to KK^{*+n}(E \te \dd A, D)
\end{equation*}
and 
\begin{equation*} 
\delta_* \cc KK^*(E \te \dd A, D) \to KK^{*+n}(E, D \te A)
\end{equation*}
are inverse isomorphisms.
\item If $A$ has a dual, and $A'$ is KK-equivalent to A, then $A'$ has
  a dual which is KK-equivalent to the dual of $A$.
\end{enumerate}
For details and further development, see ~\cite{k-p2}. 
It is not apparent if an algebra has a dual or not. Indeed, the main
goal of this paper is to exhibit an example of a class of algebras
with specific types of duals which have a geometric and dynamical
origin.  However, one can start to build up a class of algebras which
have duals in an elementary way.  For example, if $X$ is a finite
complex, then $\dd X$ exists.  If $A \in \N$ and $K_*(A)$ is finitely
generated, then $A$ is KK-equivalent to $C(X)$ where $X$ is a finite
complex, and hence $A$ has a dual.  Moreover, Connes has shown that
$\mathcal{A}_{\theta}$ is self-dual for $\theta$ irrational.  

The largest subclass of $\N$ for which $\dd$ is involutive modulo
KK-equivalence is $\dd\N$.  This can be compared with a result of M.
Boardman, ~\cite{boardman}, which states that the largest category on
which Spanier-Whitehead duality is involutive is the homotopy category
of finite complexes.  Thus, $\dd\N$ has a formal similarity with the
homotopy category of finite complexes.  This itself does not clarify
the issue of which $\cs$-algebras should play the role of
non-commutative finite complexes, but it is suggestive.  
This will be discussed further in ~\cite{k-p2}.

One may view duality as being of even or odd type depending on whether
$\Delta$ belongs to $KK^n(A \te \dd A,\C)$ for $n$ even or odd.  We
will discuss the odd type of duality here. However, in connections to
the Novikov Conjecture, ~\cite{ kaminker-p1}, and physics, ~\cite{ connes:book}, the even
type naturally appears. 

\section{Criterion for duality classes}

In this section we will present a technical result, Proposition ~\ref{hyp}, which gives a criterion for when two classes $\Delta$ and
$\delta$ yield duality isomorphisms. This is essentially the same as
the condition given by Connes, ~\cite[p. 588]{connes:book}, except
that our duality is in the odd case and this requires adjusting the
arguments for Bott periodicity.  This technicality is actually what
allows us to obtain the duality isomorphisms in the case of shifts of
finite type.  

Thus, we shall give useable conditions under which $\ds \co \de = 1$
and $\de \co \ds = 1$.  This breaks into two parts.  The first is an
uncoupling step and the second is a type of cancellation.  In the
following sections we apply this to the case of Cuntz-Krieger
algebras.

We will first prove in detail that $\delta_0 \Delta_0 = 1_{K_0(A)}$.
The statement, $\delta_1 \Delta_1 = 1_{K_1(A)}$, follows in a similar
manner.  We then sketch the proof that $ \Delta_0 \delta_0 =
1_{K_0(B)}$.  To start with we will perform the uncoupling step.  Let $x
\in K_0(A) =  KK(\C,A)$.  Then we have
\begin{equation*} 
  \delta_0 \Delta_0 (x) = \be \te _{\sk \te \sk} (\delta \te_{B}
  (x \te_{A} \Delta)).
\end{equation*}
Consider, first, the factor $(\delta \te_{B} (x \te_{A}
\Delta))$.  We have
\begin{align*}
  (\delta \te_{B} (x \te_{A} \Delta)) &= \tk_{\sk} (\delta) \te (\tk^{A}
  \tk_{\sk} \tk_{B} (x) \te \tk^{A} (\Delta)) \\ &= (\tk_{\sk}
  (\delta) \te (\tk^{A} \tk_{\sk} \tk_{B} (x)) \te \tk^{A} (\Delta).
\end{align*}
Now, a direct computation yields that
\begin{align*}
  (\tk_{\sk} (\delta) \te (\tk^{A} \tk_{\sk} \tk_{B} (x)) \te \tk^{A}
  (\Delta) &= (\tk_{\sk} \tk^{\sk} (x) \te {\s_{12} \s^{24}} \tk^{A}
  \tk^{\sk} (\delta)) \te \tk^{A} (\Delta) \\ &= \tk_{\sk} \tk^{\sk}
  (x) \te ({\s_{12} \s^{24}} \tk^{A} \tk^{\sk} (\delta)) \te \tk^{A}
  (\Delta)).
\end{align*}
Putting $\be$ back into the product and simplifying, one obtains
\begin{equation*} 
  \be \te _{\sk \te \sk} (x \te_{A} (\delta \te_B \Delta)) = x \te_A
  (\be \te _{\sk \te \sk} (\delta \te_B \Delta)).
\end{equation*}
This accomplishes the uncoupling.

\begin{prop}
One has $\delta_0 \Delta_0 (x) = x \te_A (\be \te _{\sk \te \sk}
(\delta \te_B \Delta))$.
\end{prop}
What one would hope is that
\begin{equation} 
\label{hope}
\be \te _{\sk \te \sk} (\delta \te_B \Delta) = 1_A \in KK(A,A),
\end{equation}
thus yielding
\begin{equation*}
\delta_0 \Delta_0 (x) = x.
\end{equation*}
Indeed, if $\delta \te_B \Delta = \tk_A ( \T_0)$, then $\be _{\sk \te
\sk} \tk^A (\T_0) = 1_A$ by Bott periodicity.  However, this need not be the case.  This is because $\delta$ and
$\Delta$ behave like K-theory fundamental classes and may differ by a
unit from ones which would yield ~\eqref{hope}.  There is a way to
compensate for this which we address next.

\begin{prop}
\label{hyp}
Suppose that there are automorphisms $\Theta_A \cc A \te \sk \to A \te
\sk$ and $\Theta_B \cc B \te \sk \to B \te \sk$ such that 
\begin{align*}
  (\Theta_A )_i \cc K_i(A \te \sk) \to K_i(A \te \sk) \\ (\Theta_B )_i
  \cc K_i(B \te \sk) \to K_i(B \te \sk)
\end{align*}
are the identity map, for $i= 0,1$, and, further,
\begin{align}
\label{thiii}
\s_{12}(\delta \te_B \s_{12}(\Delta)) &= \Theta_A \te_{A \te \sk} \tk^A
(\T_0) \\
\s_{12}(\delta \te_A \s_{12}(\Delta)) &= \Theta_B \te_{B \te \sk} \tk^B
(\T_0) 
\end{align}
Then,
\begin{equation*} 
\delta_i \Delta_i \cc K_i(A) \to K_i(A)
\end{equation*}
is the identity for $i=0,1$
\end{prop}
\begin{proof}
We will give the proof for $\delta_0 \Delta_0$, the other case being
similar.
Condition ~\eqref{thiii} states that 
\begin{equation*} 
\s_{12} \tau_{\sk} \tau_{A} (\delta) \te (\tau_{A} (\s_{12}(\Delta))) =
\tau^{\sk} (\Theta_A ) \te \tau^{A} (\T_0).
\end{equation*}
Thus, one has
\begin{equation*} 
\be \te _{\sk \te \sk} (\delta \te_B (\s_{12}(\Delta))) = \tk^A (\be) \te
\tk_{\sk}(\Theta_A) \te \tk^A (\T_0). 
\end{equation*}
Now, 
\begin{align*}
\tk^A(\be) \te \tk_{\sk} (\Theta_A) &= (\Theta_A)_* (\tk^A(\be)) \\
&= \tk^A (\be),
\end{align*}
so we obtain
\begin{align*}
\be \te _{\sk \te \sk} (\delta \te_B \s_{12}(\Delta) &= \tk^A (\be) \te
\tk^A (\T_0) \\
&= 1_A,
\end{align*}
which yields the desired result.
\end{proof}

For the composition $\Delta_* \delta_*$ we have a similar result.
\begin{prop}
Under the hypothesis of Proposition ~\ref{hyp}, we have that 
\begin{equation} 
\Delta_i \delta_i \cc K^{i+1}(B) \to K^{i+1}(B)
\end{equation}
is the identity for $i=0,1$.
\end{prop}
\begin{proof}
The proof is obtained from the previous one by making obvious changes.
\end{proof}

The other cases follow in the same way.  Thus,  showing that one has a
duality between algebras reduces to constructing the maps $\Theta_A$
and $\Theta_B $ satisfying the conditions above.  In the next two sections
we will do this for the case of the stable and unstable Ruelle
algebras associated to a subshift of finite type.     

\section{Construction of duality classes for shifts of finite type}

In this section we will construct the classes in KK-theory needed to
exhibit the duality between $\oa \te \K$ and $\oat \te \K$.  Let $A$
be an $n \x n$ matrix with entries which are all zero or one.  We
assume that $A$ has no row or column consisting entirely of
zeros and that the associated shift space is a Cantor set.  

The Cuntz-Krieger algebra, $O_A$, is the universal $\cs$-algebra
generated by partial isometries $s_1, \dots , s_n$ satisfying 
\begin{enumerate}[i)]
\item the projections $s_1 s_{1}^{*}, \dots , s_n s_{n}^{*}$ are
pairwise orthogonal and add up to the identity of $O_A$,
\item for $k = 1, \dots ,n$ one has
\begin{equation} 
s_{k}^{*} s_k = \sum_i A_{ki} s_i s_{i}^{*}.
\end{equation}
\end{enumerate}
The condition above, that the shift space be a Cantor set, guarantees
that the algebra described does not depend on the choice of the
partial isometries, ~\cite{cuntz-k1}.
If $A_{ij} = 1 $ for all $i,j$, then the algebra $O_A$ is denoted
$O_n$.

In a similar manner we consider $O_{A^T}$, with generators $t_1,
\dots, t_n$ satisfying
\begin{equation} 
t_{k}^{*} t_k = \sum_i A_{ik} t_i t_{i}^{*}
\end{equation}
for $k = 1, \dots ,n$.

Our aim in this section is to explicitely construct the elements
\begin{equation*}
\delta \in KK^1(\C, O_A \te O_{A^T})
\end{equation*}
and 
\begin{equation*}  
\Delta \in KK^1(O_A \te O_{A^T}, \C)
\end{equation*}
which are needed to show that $O_A$ and $O_{A^T}$ are dual.

The construction of $\delta$ is the easier of the two, (c.f. ~\cite{cuntz-k2}).
Let
\begin{equation} 
w = \sum_{i=1}^{n} s_{i}^{*} \te t_i \in \oa \te \oat.
\end{equation}
Then one has
\begin{equation*} 
w^* w = w w^* = \sum_{i,j} A_{ij} s_j s_{j}^{*} \te t_i t_{i}^{*}.
\end{equation*}
We let $\bar w \cc C(S^1) \to \oa \te \oat$ denote both the
(non-unital) map defined by
\begin{equation} 
\bar w (z) = w
\end{equation}
as well as its restriction to $C_0(\R) \subseteq C(S^1)$.  
\begin{definition}
Let $\delta \in KK^1(\C, \oa \te \oat)$ be the element determined by
the homomorphism $\bar w$.
\end{definition}
The element $\Delta$ is constructed using the full Fock space of a finite
dimensional Hilbert space.  (For related constructions see the papers
of D. Evans and D. Voiculescu, ~\cite{voiculescu,evans}.)

Let $\h$ denote an n-dimensional Hilbert space with orthonormal basis
$\xi_1, \dots, \xi_n$.  Let $\h^{\te m} = \h \te \cdots \te \h$ be the
m-fold tensor product of $\h$ and let $\h_0$ be a one dimensional
Hilbert space with unit vector $\Omega$.  Then the full Fock space of $\h$,
$\F$, is defined to be 
\begin{equation*} 
\F = \h_0 \oplus (\bigoplus_{n=1}^{\infty} \h^{\te n})
\end{equation*}
There is a natural orthonormal basis for $\F$,
\begin{equation*} 
\{\Omega, \xi_{i_1}\te \cdots \te \xi_{i_m}|m=1,2,\dots,\qquad  1\leq i_j
\leq n\}. 
\end{equation*}
Define the left and right creation operators, $L_1, \dots , L_n$ and
$R_1, \dots ,R_n$, on $\F$ by 
\begin{equation*} 
L_k \Omega = \xi_k = R_k \Omega
\end{equation*}
and
\begin{align}
L_k(\xi_{i_1}\te \cdots \te \xi_{i_m} ) &= \xi_k \te \xi_{i_1}\te
\cdots \te \xi_{i_m} \\
R_k(\xi_{i_1}\te \cdots \te \xi_{i_m} ) &= \xi_{i_1}\te \cdots \te
\xi_{i_m} \te \xi_k
\end{align}
Next, we bring in the matrix $A$.  Let $\F_A \subseteq \F$ denote the
closed linear span of the vectors $\Omega$ and those $\xi_{i_1}\te
\cdots \te \xi_{i_m} $ satisfying the condition that
$A_{{i_j},{i_{j+1}}} = 1$ for all $j = 1, \dots , m-1$.  Let $P_A$
denote the orthogonal projection of $\F$ onto $\F_A$.  Let 
\begin{align*}
L_k^A &= P_A L_k P_A \in \B(\F_A) \\
R_k^A &= P_A R_k P_A \in \B(\F_A)
\end{align*}
for $k= 1, \dots, n$.

It is easily checked that one has the following formulas.
\begin{align*}
L_k^A \xi_{i_1}\te \cdots \te \xi_{i_m} &= A_{k,i_1} \xi_k \te \xi_{i_1}\te
\cdots \te \xi_{i_m} \\
R_k^A \xi_{i_1}\te \cdots \te \xi_{i_m} &= A_{i_m,k}  \xi_{i_1}\te
\cdots \te \xi_{i_m} \te \xi_k \\
(L_k^A)^* \xi_{i_1}\te \cdots \te \xi_{i_m} &= A_{k,i_1}  \xi_{i_2}\te
\cdots \te \xi_{i_m} \\ 
(R_k^A)^* \xi_{i_1}\te \cdots \te \xi_{i_m} &= A_{i_m,k}  \xi_{i_1}\te
\cdots \te \xi_{i_{m-1}}. 
\end{align*}
From this one easily obtains the following result.
\begin{prop}
\label{formulas}
The operators $R^A_k$ and $L^A_k$ are partial isometries and satisfy
\begin{enumerate}[i)]
\item $(L^A_k)^* L^A_k = \sum_i A_{ki} L^A_i (L^A_i)^* + P_\Omega$
\item $(R^A_k)^* R^A_k = \sum_i A_{ik} R^A_i (R^A_i)^* + P_\Omega$
\item $[L^A_k, R^A_l] = 0$
\item $[(L^A_k)^* , R^A_l] = \delta_{kl} P_\Omega$
\end{enumerate}
\end{prop}
We are now able to construct the element $\Delta$.  Let $\E \subseteq
\B(\F_A)$ be the $\cs$-algebra generated by $\{R^A_1, \dots, R^A_n,
L^A_1, \dots , L^A_n\}$.  By Proposition~\ref{formulas} the operator
$P_\Omega$, which is compact, is in $\E$.  It is easy to check that
there is no non-trivial $\E$-invariant subspace of $\F_A$.  Thus, $\E$
contains the compact operators, $\K(\F_A)$.

Modulo the ideal $\K(\F_A)$ the elements  $L^A_1, \dots , L^A_n$ and
$R^A_1, \dots, R^A_n$ satisfy the relations for $\oa$ and $\oat$
respectively.  Moreover, the $L^A_i$'s and the $R^A_j$'s commute
modulo $\K(\F_A)$.  It follows that the $\cs$-algebra $\E / \K(\F_A)$
is a quotient of $\oa \te \oat$.  In fact they are isomorphic.  This
follows since both $\oa$ and $\oat$ are nuclear and the ideal
structure of their tensor product may be completely described in terms
of the ideals of $\oa$ and $\oat$.  These, in turn, have been
completely described in ~\cite{ cuntz-k2}.  It is then straightforward to
verify that the generators of the ideals of $\oa \te \oat$ give rise
to non-compact operators (via the $L^A_k$ and $R^A_k$) and thus $\E
/\K(\F_A) \iso \oa \te \oat$.
\begin{definition}
Let $\Delta \in KK^1(\oa \te \oat, \C)$ be the class determined by the
exact sequence
\begin{equation} 
\begin{CD}
0 @>>> \K(\F_A) @>>> \E @>{\pi_{A}}>> \to \oa \te \oat \to 0.
\end{CD}
\end{equation}
\end{definition}
Note that one has
\begin{equation*} 
\pi_A (R^A_k) = 1 \te t_k
\end{equation*}
and 
\begin{equation*} 
\pi_A (L^A_k) = s_k \te 1.
\end{equation*}


 \section{Duality for Cuntz-Krieger algebras}
In this section we will show that the duality classes constructed in
the previous section actually implement a duality isomorphism for the
algebras $\oa$ and $\oat$.  According to Proposition ~\ref{hyp}, it
will be sufficient to construct homomorphisms
\begin{equation}
\begin{align}
\Theta_{\oa} \cc \oa \te \sk \to \oa \te \sk \\
\Theta_{\oat} \cc \oat \te \sk \to \oat \te \sk
\end{align}
\end{equation}
which satisfy the conditions stated there.  That is, we
must show that $\Theta_{\oa}$ and $\Theta_{\oat}$ induce the identity
homomorphism on K-theory and satisfy the second condition in
Proposition ~\ref{hyp} which states
\begin{align*}
\s_{12}(\delta \te_{\oat} \s_{12}(\Delta)) &= \Theta_{\oa} \te_{{\oa} \te \sk} \tk^{\oa}
(\T_0) \\
\s_{12}(\delta \te_{\oa} \s_{12}(\Delta)) &= \Theta_{\oat} \te_{{\oat} \te \sk} \tk^{\oat}
(\T_0)
\end{align*}
We will work out the details only for $\Theta_{\oa}$, the other case
being similar.

To define $\Theta_{\oa}$ we first set
\begin{equation*}
\bar{\Theta} \cc \oa \te C(S^1) \to \oa \te C(S^1)
\end{equation*}
by
\begin{align*}
\bar{\Theta} (1 \te z) = 1 \te z \\
\bar{\Theta} (s_i \te 1) = s_i \te z.
\end{align*}
Then $\bar{\Theta} $ extends to an automorphism of $\oa \te
C(S^1)$, as follows from the universal property of $\oa$.
The diagram
\begin{equation*}
\begin{CD}
\oa \te C(S^1) @>{\bar{\Theta}}>> \oa \te C(S^1)\\
@V{1_{\oa}} \te \pi VV    @V{1_{\oa}} \te \pi VV \\
\oa @>id>> \oa
\end{CD}
\end{equation*}
commutes, where $\pi \cc C(S^1) \to \C$ is defined by $\pi_1 (z) =
1$.  It follows that we may define $\Theta_{\oa} = \bar{\Theta} |
\ker(1_{\oa} \te \pi)$.  It is
an automorphism of $\oa \te \sk$.  We now must show that
$\Theta_{\oa}$ satisfies the necessary conditions.
\begin{theorem}
The maps
\begin{equation*}
{\Theta_{\oa}}_* \cc K_i(\oa \te \sk) \to K_i(\oa \te \sk)
\end{equation*}
are the identity for $i = 0 , 1$
\end{theorem}
\begin{proof}
Recall that
\begin{equation*}
\oa \te \K \iso \bar F_A \rtimes_{\sigma_A} \Z
\end{equation*}
where $\bar F_A$ is a stable $AF$-algebra with automorphism
$\sigma_A$.  In this situation, $\oa$ is actually a full corner in
$\bar \oa$ and compressing $\bar F_A$ to this corner yields $\bar F_A \subseteq \oa $ which is the
closure of the ``balanced words'' in the $s_i$'s as described in
~\cite{ cuntz-k1}.  Observe that the restriction of $\bar{\Theta} $ to
$F_A \te C(S^1)$ is the identity.  We will apply the
Pimsner-Voiculescu exact sequence to compute $K_*(\oa \te \sk)$,
making necessary modifications since $\bar F_A$ is not unital and then
study ${\Theta_{\oa}}_*$.

Let $B$ denote the multiplier algebra of $\bar F_A \te \sk \te \K$
where $\K = \K(l^2(\NN))$.  Let $e_{ij}$ denote the standard matrix units in
$\K$.  Define $\rho \cc \bar F_A \te \sk \to B$ by
\begin{equation*}
\rho ( a \te b) = \sum_{i \in \NN} \s^{i}_A (a) \te f \te e_{ii}
\end{equation*}
where the sum is taken in the strict topology.  Let $S$ denote the
unilateral shift on $\ell^2(\NN)$.  Let $D$ denote the
$\cs$-algebra generated by $\bar F_A \te \sk \te \K$, $1 \te 1 \te S$
and $\{ \rho(a \te f) | f \in \sk, a \in \bar F_A\}$.  Let $D_0$ be
the ideal in $D$ generated by $\bar F_A \te \sk \te \K$ and $\{ \rho(a
\te f) | f \in \sk, a \in \bar F_A\}$.  There is an exact sequence
\begin{equation*}
0 \to \bar F_A \te \sk \te \K \to D_0 \to \sk \te (\bar F_A \cp \Z)
\to 0.
\end{equation*}
Moreover, the two maps
\begin{equation*}
j \cc \bar F_A \te \sk \to \bar F_A \te \sk \te \K
\end{equation*}
defined by $j(a \te f) = a \te f \te e_{11}$
and
\begin{equation*}
\rho \cc \bar F_A \te \sk \to D
\end{equation*}
both induce isomorphisms on K-theory.

Finally, we have
\begin{equation}
K_0(\bar F_A \te \sk) \iso K_1(\bar F_A)= 0
\end{equation}
since $\bar F_A$ is an AF-algebra.  Putting this together, we obtain
the Pimsner-Voiculescu sequence for the $\oab$'s:
\begin{equation*}
0 \to K_1(\oab \te \sk) \to K_1(\bar F_A \te \sk) \to K_1(\bar F_A \te
\sk) \to K_0(\oab \te \sk) \to 0
\end{equation*}

We define an automorphism $\tilde \Theta$ of $D$ by
\begin{equation}
\tilde \Theta = ad (\sum_{i \in \NN} 1 \te z^i \te e_{ii})
\end{equation}
where, again, the sum is in the strict topology.  Notice that $\tilde
\Theta \co \rho = \rho$ and $\tilde \Theta | (\bar F_A \te \sk \te \K)$
is approximately inner and hence trivial on K-theory.  Also observe
that $\tilde \Theta | (\bar F_A \te \sk \te \K) = \bar F_A \te \sk \te
\K$ and that the automorphism of the quotient of $D_0$ by $\bar F_A
\te \sk \te \K$ induced by $\tilde \Theta $ is precisely $\Theta_{\oa}
$,
after identifying this quotient with $\oab \te \sk$ and restricting to
$\oa \te \sk \subseteq \oab \te \sk$.  We have a commutative diagram
\begin{equation*}
\begin{CD}
0 @>>> K_0(\oab \te \sk) @>>> K_1(\bar F_A \te \sk) @>>> K_1(\bar F_A
\te \sk) @>>> K_1(\oab \te \sk) @>>>0 \\
@. @VV{{\Theta_{\oa}}_*}V   @VV{\tilde \Theta_*}V  @VV{\tilde \Theta_*}V
@VV{{\Theta_{\oa}}_*}V\\
0 @>>> K_0(\oab \te \sk) @>>> K_1(\bar F_A \te \sk) @>>> K_1(\bar F_A
\te \sk) @>>> K_1(\oab \te \sk) @>>>0
\end{CD}
\end{equation*}
From the observations above, we have both maps $\tilde \Theta_* = id$,
and it follows that ${\Theta_{\oa}}_* $ is the identity.
\end{proof}
It remains for us to verify that  condition
~\eqref{thiii} is satisfied.  To that end we observe first that
\begin{equation*}
\s_{12}(\delta \te_{\oat} \s_{12}(\Delta)) = \Theta_{\oa} \te_{{\oa} \te \sk} \tk^{\oa}
(\T_0)
\end{equation*}
is equivalent to
\begin{equation}
\label{81}
\tk_{\sk}((\Theta_{\oa})^{-1}) \te \s_{12} \tk_{\sk} \tk_{\oa}
(\delta) \te \tk^{\oa} (\s_{12}(\Delta) ) =  \tk^{\oa}
(\T_0).
\end{equation}
Thus, we will prove the latter statement.

Now, $\tau^{\oa}(\s_{12}(\Delta) ) \in KK^1(\oa \te \oat \te \oa, \oa)$ was obtained from the
extension
\begin{equation*}
0 \rar{}{} \K \te \oa \rar{}{} \E \te \oa \rar{\pi_A \te 1_{\oa}}{} \oa \te \oat \te \oa \rar{}{} 0.
\end{equation*}
Moreover, the remaining term
\begin{equation*}
\tk_{\sk}((\Theta_{\oa})^{-1}) \te \s_{12} \tk_{\sk} \tk_{\oa} (\delta)
\end{equation*}
actually yields a $*$-homomorphism from $ \oa \te \sk \te \sk$ to $
\oa \te \oat \te \oa \te \sk$.  Thus, the left side of ~\eqref{81} is
represented by applying $\tk_{\sk}$ to the element represented by the top row of the following diagram
\begin{equation*}
\begin{CD}
0 @>>> \K \te \oa @>>> \E' @>>> \oa \te \sk @>>>0 \\
@. @VVV @VVV @VV{1_{\oa} \te i}V \\
0 @>>> \K \te \oa @>>> \E'' @>>> \oa \te C(S^1) @>>>0 \\
@. @VVV @VVV @VV{\bar \al}V \\
0 @>>> \K \te \oa @>>> \E \te \oa @>>> \oa \te \oat \te \oa @>>> 0
\end{CD}
\end{equation*}
where $\al = (\Theta_{\oa})^{-1} \te  \tk_{\oa} (\delta) = \bar
\al \co (1_{\oa} \te i)$, and $\al \te 1_{\sk} = \tk_{\sk}((\Theta_{\oa})^{-1}) \te \s_{12} \tk_{\sk} \tk_{\oa} (\delta)$.

The crucial step is to untwist the middle row by finding an isomorphism $\E'' \iso \T \te \oa$ so
that the following diagram commutes
\begin{equation*}
\begin{CD}
0 @>>> \K \te \oa @>>> \E'' @>>> \oa \te C(S^1) @>>> 0 \\
@. @V{=}VV @V{\iso}VV @V{\s_{12}}VV \\
0 @>>> \K \te \oa @>>> \T \te \oa @>>{\pi_{\T}} \te 1_{\oa}> C(S^1) \te \oa @>>> 0,
\end{CD}
\end{equation*}
where $\T$ is the Toeplitz extension.

Assuming this, the proof can be completed as follows.  We have
\begin{equation*}
\tk_{\oa}(\T) = \s_{12} \bar \alpha^* (\tk^{\oa} (\s_{12}(\Delta) )).
\end{equation*}
  Hence, one has
\begin{align*}
\tk_{\oa}(\T_0) &=  (i \te 1_{\oa})^* (\tk^{\oa}(\T)) \\
&= (i \te 1_{\oa})^* \s_{12} \bar \alpha^* (\tk^{\oa}(\s_{12}(\Delta) )) \\
&= \s_{12} (1_{\oa} \te i)^* \bar \alpha^* (\tk^{\oa}(\s_{12}(\Delta) )) \\
&= \s_{12} \alpha^* (\tk^{\oa}(\s_{12}(\Delta) )).
\end{align*}
Thus, substituting in for $\al$, we obtain
\begin{equation*}
\tk_{\oa}(\T_0) = \tk_{\sk} ((\Theta_{\oa})^{-1}) \te \s_{12}
\tk_{\sk} \tk^{\oa}(\delta) \te \tk^{\oa}(\s_{12}(\Delta) ),
\end{equation*}
which is the desired formula.

We now turn to the issue of obtaining the explicit isomorphism between
$\E''$ and $\T \te \oa$.  For convenience, we will suppress the $A$ in our
notation from the elements such as ${R_i}^A$, and ${L_i}^A$.  Define $W$ in $\E \te
\oa$ by
\begin{equation*}
W = \sum_{i=1}^{n} R_i \te {s_i}^*.
\end{equation*}
We will need two technical lemmas.
\begin{lem}
\label{43}
One has
\begin{enumerate}[i)]
\item $\pi \te 1_{\oa} (W) = \bar \al (1 \te z)$.
\label{i}
\item $W^* W = \sum_{i,j} A_{ji} R_j {R_j}^* \te s_i {s_i}^* +
P_{\Omega} \te 1$.
\item $[W^* , W] = P_{\Omega} \te 1$.
\item $(P_{\Omega} \te 1) W = 0$.
\item $[W, L_k \te 1] = 0$ for $k = 1,\dots,n$.
\item $[W^* , L_k \te 1] = P_{\Omega} \te s_k$ for $k = 1,\dots,n$.
\end{enumerate}
\end{lem}
\begin{proof}
For (\ref{i}), one proceeds as follows.  Note first that 
\begin{align*}
  \bar \alpha ( 1 \te z) &= \s^{23} (\bar \Theta_{\oa}^{-1})^* (1
  \te \bar \omega (z)) \\
&= \s^{23} (1 \te \bar \omega (z)) \\
&= \s^{23} ( \sum_i 1 \te s_{i}^{*} \te t_i) \\
&= \sum_i 1 \te t_i \te s_{i}^{*}.
\end{align*}
Moreover, 
  \begin{align*}
    (\pi \te 1_{\oa})(W) &= (\pi \te 1_{\oa}) (\sum_i R_i \te
    s_{i}^{*}) \\
    &= \sum_i \pi(R_i) \te s_{i}^{*} \\
    &= \sum_i 1 \te t_i \te s_{i}^{*}.
  \end{align*}
The remaining parts of lemma can be verified in a routine manner.
\end{proof}

The remaining facts we need are incorporated into the following.
\begin{lem}
\label{44}
Let $V_k = W^* (L_k \te 1)$, for $k = 1 \dots n$.  Then we have, for
each $k$,
\begin{description}
\item[i)] $\pi \te 1_{\oa} (V_k) = \bar \al (s_k \te 1)$,
\label{ii}
\item[ii)] $\sum_j V_j V^{*}_{j} = W^* W$,
\item[iii)] $V_{k}^{*} V_k = \sum_j A_{kj} V_j V^{*}_{j} $,
\item[iv)] $[W, V_k] = 0$,
\item[v)] $[W^* , V_k] = 0$.
\end{description}
\end{lem}
\begin{proof}
  As in the previous lemma, we will verify (\ref{ii}) and leave the
  remaining parts of the proof to the reader, since they are
  essentially routine.  For ~\ref{ii}, we check
    \begin{align*}
  (\pi \te 1_{\oa}) (V_k) &= (\pi \te 1_{\oa})(W^*)(\pi \te
  1_{\oa})(L_k \te 1) \\
  &= \bar \alpha (1 \te z) (s_k \te 1 \te 1),
  \end{align*}
  and  
      \begin{align*}
      \bar \alpha (s_k \te 1) &= \s^{23} (1_{\oa} \te \bar w)
      (\bar \Theta_{A}^{-1})^* (s_k \te 1)  \\
      &= \s^{23} (1_{\oa} \te \bar w) (s_k \te 1)  \\
      &= \s^{23} (1_{\oa} \te w) (s_k \te 1 \te 1)  \\
      &= \bar \alpha (1 \te z) (s_k \te 1 \te 1).
    \end{align*}
  \end{proof}

Now we may define the isomorphism from $\T \te \oa$ to $\E ''$.  Let
$S$ denote the unilateral shift. The required map is defined by
sending $S \te 1$ to $W$, and $1
\te S_k$ to $V_k$.  Note that the unit of $\T \te \oa$ is mapped to
$W^* W$ in $\E''$.  The fact that this assignment extends to a *-homomorphism
follows from the universal properties of $\T$ and $\oa$.  The fact
that it is onto follows from observing that $\E ''$ is generated by $\{
W, V_1, \dots , V_k \}$ which is straightforward.  Finally, the
fact that the appropriate diagram commutes follows from ~\ref{44} and
~\ref{43}.

\section{Final comments}
\begin{enumerate}
\item As mentioned earlier, the duality theorem holds for the stable
  and unstable Ruelle algebras associated to a Smale
  space,~\cite{kaminker-p2}.  However, the duality classes, $\Delta$
  and $\delta$, must be constructed in a different way.  This is done
  using asymptotic morphisms and uses the fact that locally the Smale
  space decomposes into a product of expanding and contracting sets.
  It would be very interesting to have a Fock space construction of
  the more general classes as well.
\item The duality result for Cuntz-Krieger algebras sheds some light
  on the computations of the K-theory of $\oa$'s as in
  ~\cite{cuntz-k2}.  Recall that if $A$ is an $n \x n$ aperiodic
  matrix of $0$'s and $1$'s , then there are {\em canonical}
  isomorphisms
\begin{align*}
K_0(\oa) \iso \Z^n / (1-A^T) \Z^n \\
K_1(\oa) \iso \ker(1-A^T) \\
K^0(\oa) \iso \ker(1-A) \\
K^1(\oa) \iso \Z^n / (1-A) \Z^n.
\end{align*}
Note that $\Z^n / (1-A) \Z^n \iso \Z^n / (1-A^T) \Z^n$ by the
structure theorem for finitely generated abelian groups, but the
isomorphism is not natural.
The explanation for why one has $A^T$ in the formulas now comes from
duality, since one has the diagram
\begin{equation*}
\begin{CD}
K_0(\oa) @>\iso>> K^1(\oat) \\
@VVV  @VVV \\
\Z^n / (1-A^T) \Z^n @>{=}>> \Z^n / (1-A^T) \Z^n.
\end{CD}
\end{equation*}
\end{enumerate}

\end{document}